\documentclass[12pt]{amsart}

\usepackage{amsmath, amsthm, amssymb}
\usepackage{amssymb}
\usepackage{graphicx}
\usepackage{color}
\usepackage{transparent}
\usepackage{float}
\usepackage{wasysym}
\usepackage{subfigure}

\newtheorem{thm}{Theorem}

\newtheorem*{question}{Question}

\newtheorem*{theononumber}{Theorem}

\newcommand{\R}{\mathbf{R}}

\begin{document}
\title{Note on the distortion of $(2,q)$-torus knots}
\author{Luca Studer}
\email{luca.studer@math.unibe.ch}
\pagestyle{plain} 
\begin{abstract} We show that the distortion of the $(2,q)$-torus knot is not bounded linearly from below. 
\end{abstract} 
\maketitle

\section{Introduction}
The notion of distortion was introduced by Gromov~\cite{GLP}. If $\gamma$ is a rectifiable simple closed curve in $\R^3$, then its distortion $\delta$ is defined as

\begin{align*}
\delta(\gamma) = \sup_{v,w \in \gamma} \frac{d_{\gamma} (v,w)}{|v-w|},
\end{align*}
where $d_{\gamma}(v,w)$ denotes the length of the shorter arc connecting $v$ and $w$ in $\gamma$ and $| \cdot |$ denotes the euclidean norm on $\R^3$. For a knot $K$, its distortion $\delta(K)$ is defined as the infimum of $\delta(\gamma)$ over all rectifiable curves $\gamma$ in the isotopy class $K$. Gromov~\cite{Gr} asked in 1983 if every knot $K$ has distortion $\delta (K) \leq 100$. The question was open for almost three decades until Pardon gave a negative answer. His work~\cite{Pa} presents a lower bound for the distortion of simple closed curves on closed PL embedded surfaces with positive genus. Pardon showed that the minimal intersection number of such a curve with essential discs of the corresponding surface bounds the distortion of the curve from below. In particular for the $(p,q)$-torus knot he obtained the following bound.

\begin{theononumber}[\cite{Pa}]
Let $T_{p,q}$ denote the $(p,q)$-torus knot. Then 
\begin{align*}
\delta(T_{p,q}) \geq \frac{1}{160}  \min(p,q).
\end{align*}
\end{theononumber}
By considering a standard embedding of $T_{p,p+1}$ on a torus of revolution one obtains $\delta(T_{p,p+1})\leq const \cdot p$, hence for $q=p+1$ Pardons result is sharp up to constants.

An alternative proof for the existence of families with unbounded distortion was given by Gromov and Guth~\cite{GG}. In both works the answer of Gromovs question was obtained by an estimate of the conformal length, which is up to a constant a lower bound for the distortion of rectifiable closed curves. However the conformal length is in general not a good estimate for the distortion. For example one finds easily an embedding of the $(2,q)$-torus knot with conformal length $\leq 100$ and distortion $\geq q$ by looking at standard embeddings on a torus of revolution with suitable dimensions. In particular neither Pardon's nor Gromov and Guth's arguments yield lower bounds for $\delta(T_{2,q})$. While Pardon writes that surely $\lim_{q \to \infty} \delta(T_{2,q})= \infty$ and that there are to his knowledge no known embeddings of $T_{2,q}$ with sublinear distortion~\cite{Pa} [p.2], Gromov and Guth~\cite{GG} write that the distortion of $T_{2,q}$ appears to be $q$ up to constants [p.33]. In this article we show that the growth rate of $\delta(T_{2,q})$ is in fact sublinear in $q$.

\begin{thm}
\label{T1}
Let $q \geq 50$. Then $\delta(T_{2,q}) \leq 7q/\log q$. In particular the distortion of the $(2,q)$-torus knot is not bounded linearly from below.
\end{thm}

\section{acknowledgments}
Most thank is owed to Sebastian Baader and Pierre Dehornoy for their inspiring introduction to the topic. I would like to thank John Pardon for useful comments on the presented example and Peter Feller, Filip Misev, Johannes Josi and Livio Liechti for many mathematical discussions. I also thank Paul Frischknecht for the pictures.

\section{Proof}
In order to prove Theorem~\ref{T1} we need to give for every odd integer $q\geq 50$ an embedding $\gamma$ of the $(2,q)$-torus knot with distortion smaller or equal to $7q/\log q$. The idea is to use a logarithmic spiral. Let $S$ be a logarithmic spiral of unit length starting at its center $0\in \R^3$ and ending at some $u \in \R^3$. An elementary calculation shows that its distortion is equal to $1/|u|$. For another path $\alpha \subset \R^3$ of unit length and diameter $\leq 2 |u|$ with endpoints $\{v,w\} = \partial \alpha$ we get
\begin{align*}
\delta(\alpha)\geq \frac{d_{\alpha}(v, w)}{|v-w|}=\frac{1}{|v-w|} \geq \frac{1}{2|u|}=\frac{\delta(S)}{2}.
\end{align*}
Hence up to at most a factor $2$ the logarithmic spiral has the smallest distortion among all paths for a prescribed pathlength-pathdiameter-ratio. It seems therefore natural to pack the $q$ windings of the $(2,q)$-torus knot into a logarithmic spiral in order to minimize distortion.

\begin{proof} [Proof of Theorem 1] Let $q$ be an odd integer greater or equal to $50$, and $k=\log(q)/2 \pi q$. We define the embedding $\gamma$ as the union of a segment of the logarithmic spiral with slope $k$, denoted by $S$, and a piecewise linear part, denoted by $L$, see Figure~\ref{f2}. The segment of the logarithmic spiral S is contained in the yellow painted $(x,z)$ plane and parametrized by 
\begin{align*}
\varphi: [0,\pi q] \to \R^2,  \  \  \  \varphi(s)= e^{ks} \cdot
\begin{pmatrix}
\cos(s) 	 		\\
\sin(s) 	 		\\
\end{pmatrix} ,
\end{align*}
see Figures~\ref{f2} and \ref{f3}. The segment of the piecewise linear part L is in the green painted $(x,y)$ plane, see Figures~\ref{f2} and \ref{f4}. Note that $$\vert \varphi(\pi q) \vert=e^{k\pi q}=\sqrt{q} \ \ \ \text{and} \ \ \ \vert \varphi(0) \vert = 1,$$ hence the lengths defining $L$ in Figure~\ref{f4} are chosen such that the union $\gamma$ of $S$ and $L$ is the simple closed curve illustrated in Figure~\ref{f2}. The linear segments $L_1$ and $L_2$ indicated in Figure~\ref{f4} are named because of their special role in the following computations. 

\begin{figure}[h]
  \def\svgwidth{300pt}
  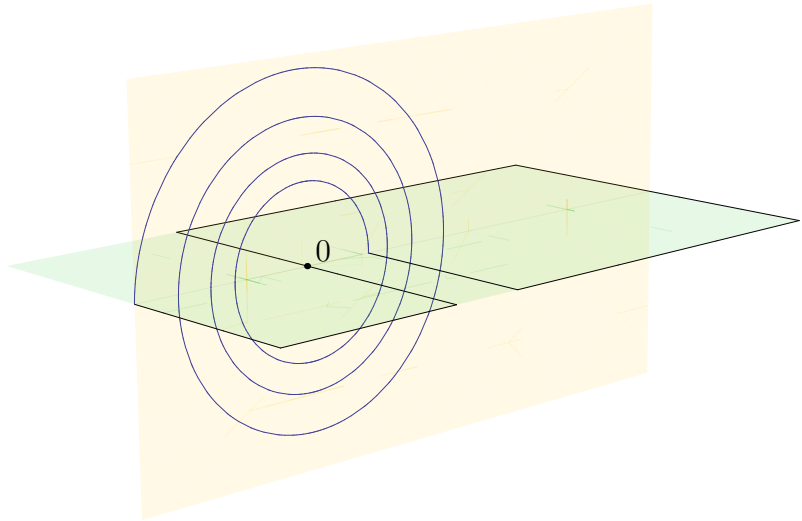
  \caption{The embedding $\gamma$ for $q=7$.}
  \label{f2}
\end{figure}

\begin{figure}[h]
\def\svgwidth{213pt}
  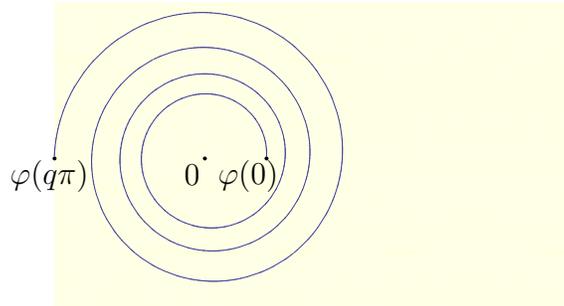
    \caption{The logarithmic spiral S in the $(x,z)$ plane.}
  \label{f3}
\end{figure}

\begin{figure}[h]
  \def\svgwidth{212pt}
  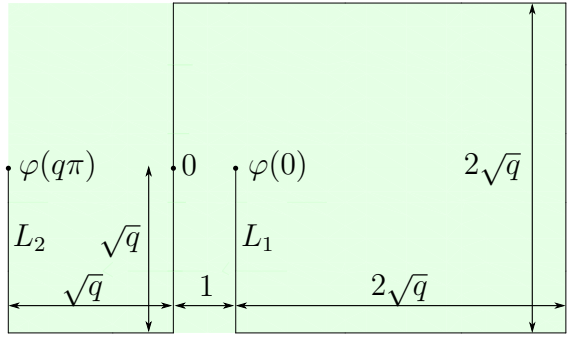
  \caption{The linear part L in the $(x,y)$ plane.}
  \label{f4}
\end{figure} \newpage \newpage \newpage

To see that the obtained curve is an embedded $(2,q)$-torus knot, we perturb $\gamma$, see Figure~\ref{f5}. This simple closed curve is ambient isotopic in $\R^3$ to $\gamma$ and if we project it onto the $(x,y)$ plane, we see a well known diagram of the $(2,q)$-torus knot, see Figure~\ref{f6}.

\begin{figure}[h]
  \def\svgwidth{300pt}
  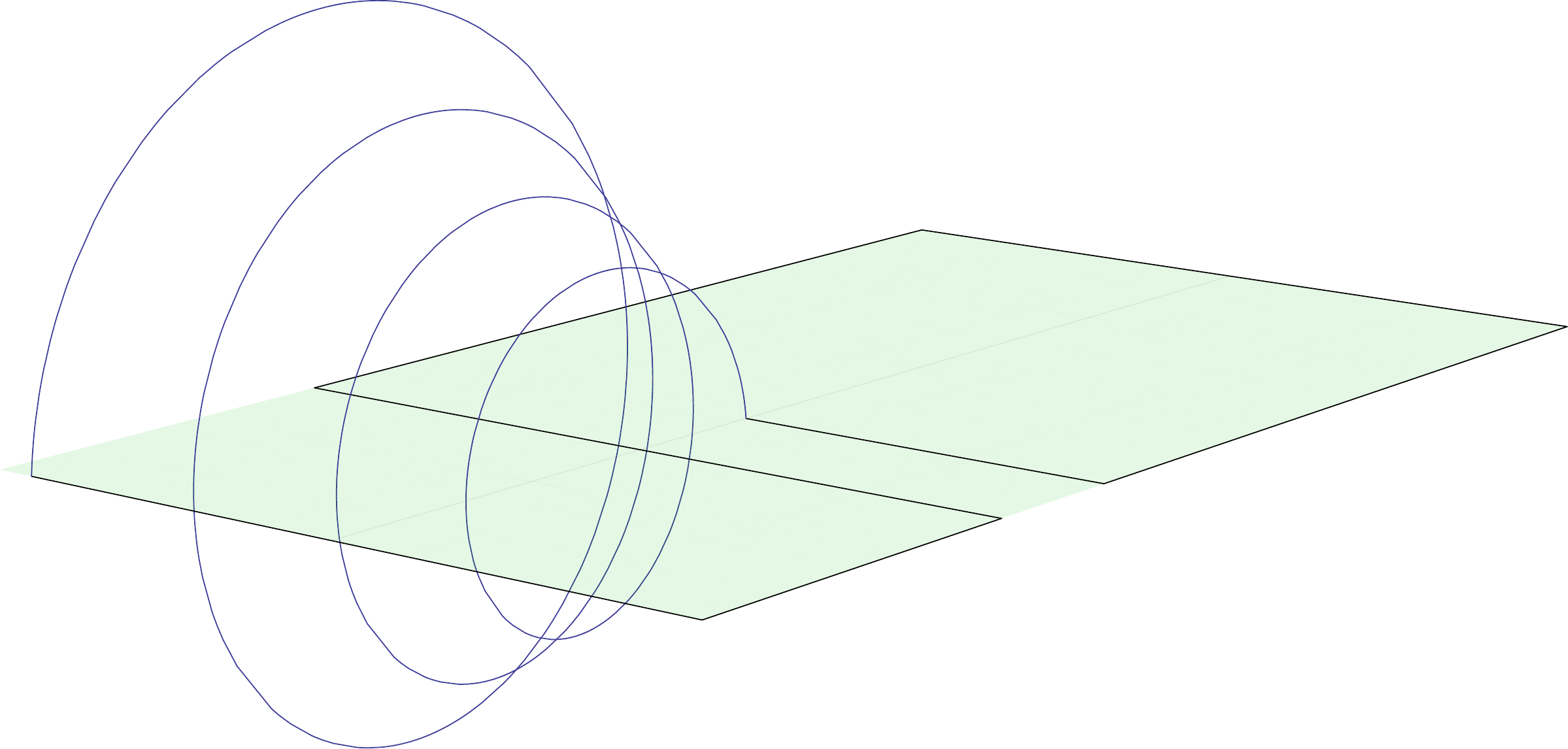
  \caption{Perturbation of $\gamma$.}
  \label{f5}
\end{figure}

\begin{figure}[h]
  \def\svgwidth{212pt}
  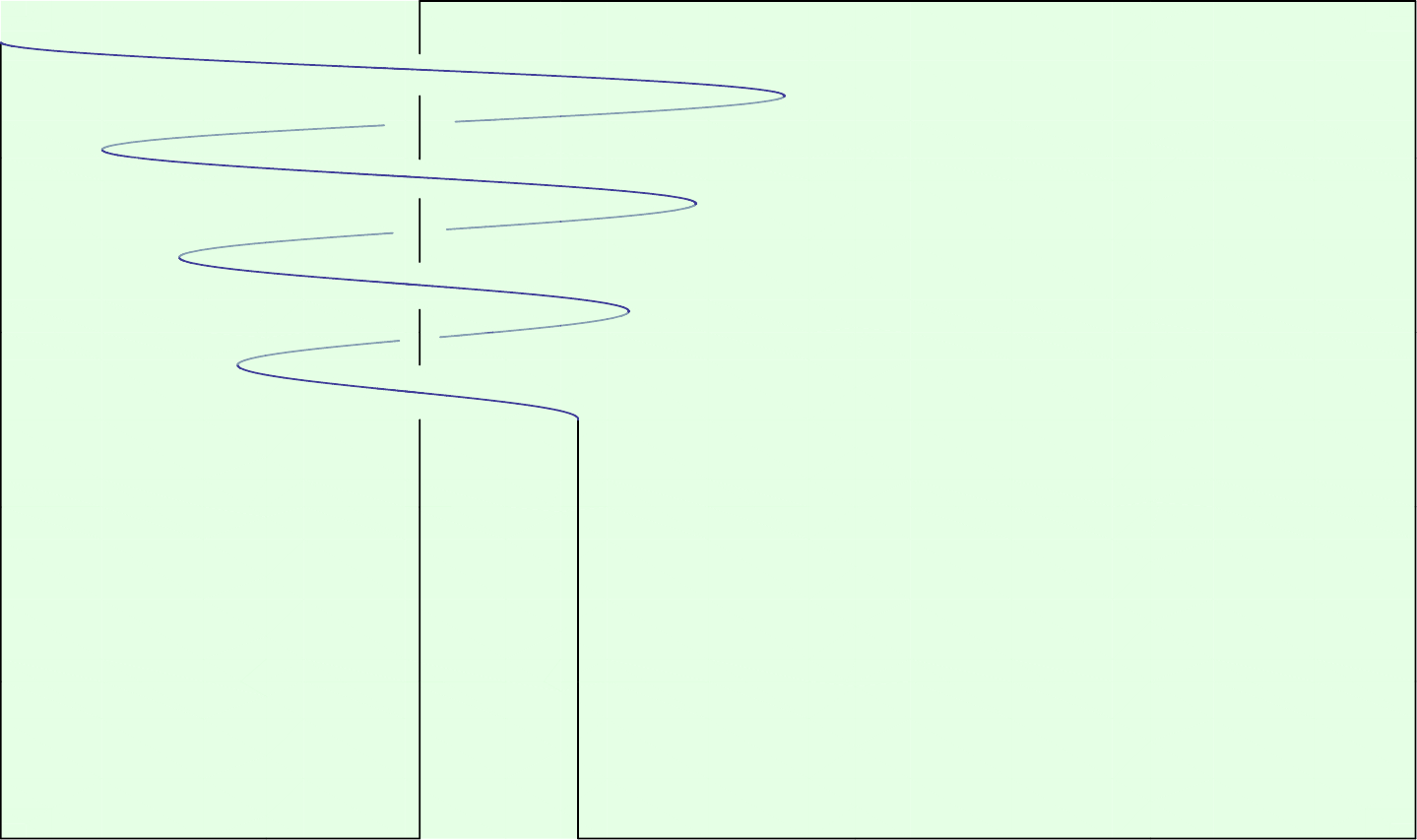
  \caption{Projection onto the $(x,y)$ plane.}
  \label{f6}
\end{figure}

We now estimate the distortion of $\gamma$. One has to show that $$\frac {d_{\gamma} (v,w)} {\vert v -w \vert}  \leq \frac {7q}{\log q}$$ for all pairs of points $v,w \in \gamma$. A calculation shows that 

$$\frac{1}{k} \cdot \sqrt{2 k^2 +1}=  \frac{2 \pi q}{\log q}   \cdot    \sqrt{2 (\log q/2 \pi q)^2 +1}  \leq \frac {7q}{\log q}$$ for all positive integers. Therefore, it suffices to show that $$\frac {d_{\gamma} (v,w)} {\vert v-w \vert}	 \leq \frac{\sqrt{2 k^2 +1}}{k}. $$ In order to do this, we distinguish four cases. \newpage			
\textit{Case 1: $v,w \in S.$} Let $0 \leq s \leq t \leq \pi q, \ v=\varphi(s), w= \varphi(t)$. From
 
\begin{align*}
\vert \varphi'(r) \vert 				= 
\left|
\begin{pmatrix}
\cos(r) 	& 	-\sin(r)	\\
\sin(r) 	& 	\cos(r)	\\
\end{pmatrix}
\begin{pmatrix}
k e^{k r} 		\\
e^{k r}  				\\
\end{pmatrix}
\right|						=
\left|
\begin{pmatrix}
k e^{k r} 		\\
e^{k r}  				\\
\end{pmatrix}
\right|						=
{\textstyle \sqrt {k^2 +1}} \cdot e^{k r} , \parbox[c]{0.2\linewidth}{}
\end{align*}
we get
\begin{eqnarray}
				d_{\gamma} (v, w)			&\leq&	d_S (v,w)   															\nonumber \\
										&=&		\int\limits_s^t \vert \varphi'(r) \vert  dr 											\nonumber\\	
										&=&		{\textstyle \sqrt{k^2 +1}}  \int\limits_s^t e^{k r} dr 								\nonumber\\	
										&=&		\tfrac{\sqrt{k^2 +1}} {k} \cdot (e^{k t} - e^{k s})									\nonumber\\
										&=&		\tfrac{\sqrt{k^2 +1}}{k} \cdot (\vert \varphi (t) \vert - \vert \varphi (s) \vert )				\nonumber\\
										&=&		\tfrac{\sqrt{k^2 +1}}{k} \cdot (\vert w  \vert - \vert v \vert ).							\nonumber
\end{eqnarray}
Since $\vert w - v \vert \geq 	\vert w \vert- \vert v \vert$,	we conclude that

\begin{align*}
\frac{d_{\gamma} (v, w)} {\vert v-w \vert}	\leq 
\frac{\sqrt{k^2 +1}}{k} \cdot \frac {(\vert w \vert - \vert v \vert )}{(\vert w \vert - \vert v \vert )} 		= 
\frac{\sqrt{k^2 +1}}{k}.
\end{align*} \newline \newline 
\textit{Case 2: $v \in L_1 \cup L_2, \ w \in S$.} We consider the case where $v\in L_1$. The idea is to find the maximum of  $$\frac {d_{\gamma} (v, w)} {\vert v -w \vert}$$ for fixed $w$ and varying $v$. Let $t = \vert v - \varphi(0) \vert$, $a = \vert \varphi(0) - w \vert$, and $b=d_S (\varphi(0),w)$, see Figure~\ref{hei}.

\begin{figure}[h]
  \def\svgwidth{360pt}
  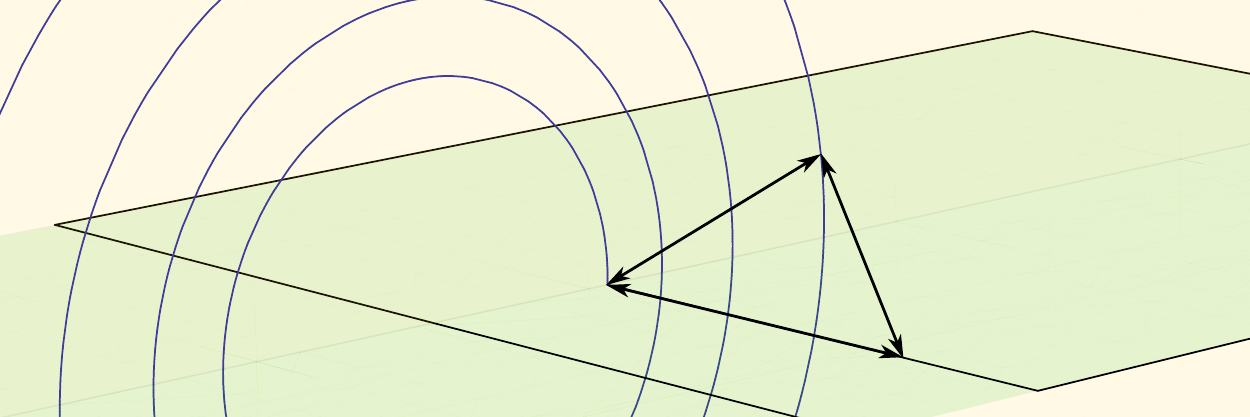
  \caption{}
  \label{hei}
\end{figure}
Note that $$\vert v - w \vert = \textstyle{\sqrt{t^2 + a^2}}$$ and $$d_{\gamma} (v,\varphi(0)) = \vert v -\varphi(0) \vert = t.$$ We get 

\begin{align*}
\frac {d_{\gamma} (v,w)} {\vert v -w \vert}			\leq 
\frac {d_{\gamma} (v,\varphi(0)) + d_S (\varphi(0),w)} {\vert v - w \vert}	= \frac {t + b} {\textstyle {\sqrt{t^2 + a^2}}}= :f(t).
\end{align*}

Deriving $f$ with respect to $t$ yields a unique critical point at $t=a^2/b$:

\begin{align*}
0=f'(t)= \frac{a^2-bt} {(a^2+t^2)^{3/2}}  \  \  \  \  \Longleftrightarrow \  \  \  \  t=a^2/b.
\end{align*}
Since $a^2/b$ is the only critical point, $f(\infty)=1 \leq b/a = f(0)$ and $$f(0) =\frac {b} {a} \leq \frac {\textstyle {\sqrt{a^2 +b^2}}} {a} =
\frac {\tfrac {a^2}{b} + b} {\textstyle{\sqrt {(\tfrac{a^2}{b})^2 + a^2}}} = f(a^2/b),$$ $a^2/b$ must be a global maximum. Consequently we get

\begin{align*}
\frac {d_{\gamma} (v, w)} {\vert v-w \vert}																&   \ \leq \			
\frac {\textstyle {\sqrt{a^2 +b^2}}} {a}  																 \\ &  \ = \ \
\sqrt{1 +  \left({\frac {b}{a}}\right)^2}										\\ &  \ = \ \
\sqrt{1 +  \left({\frac {d_S (\varphi(0),w)}{\vert \varphi(0) -w\vert}}\right)^2}										\\&\overset{\text {Case1}}{\leq}
\sqrt{1 + {\left(\tfrac{\sqrt{k^2 +1}}{k}\right)^2}}															\\&  \ = \ \
\frac{\sqrt{2 k^2 +1}}{k}.
\end{align*}

In the case where $v \in L_2$, we make the estimate with the path that connects $v$ with $w$ through $\varphi(\pi q)$. It works exactly the same and yields the same estimate. \newline \newline		
\textit{Case 3: $v,w \in L.$} Consider Figure~\ref{f4} and note that all pairs of points $v,w \in L$ that could cause big distortion are of euclidean distance at least 1. Therefore we get $$\frac{d_{\gamma} (v,w)} {\vert v-w \vert} \leq \l(L) = 11\sqrt{q}+1.$$ A calculation shows that $$11\sqrt{q} +1 \leq \frac{2\pi q}{\log q}=\frac{1}{k}$$ for q greater or equal to 50. \newline \newline 		
\textit{Case 4: $v \in L\setminus (L_1 \cup L_2), w \in S.$} Note that for these pairs of points we have $$\vert v-w \vert \geq \vert w \vert.$$ We estimate $d_{\gamma} (v,w)$ using results of Case 1 and 3:
\begin{eqnarray}
d_{\gamma} (v, w) 		&\leq&	d_{L} (v, \varphi(0)) + d_{S} (\varphi(0), w)	   														\nonumber \\
				&\leq&	 \tfrac{1} {k} +  \tfrac{\sqrt{k^2 +1}}{k} \cdot (\vert w \vert - 1) 											\nonumber\\	
				&\leq&	\tfrac{\sqrt{k^2 +1}} {k} \cdot \vert w \vert.															\nonumber	
\end{eqnarray}
We conclude that
\begin{align*}
\frac {d_{\gamma} (v, w)} {\vert v-w \vert}										\leq 
\frac {\tfrac{\sqrt{k^2 +1}} {k} \cdot \vert w \vert} {\vert w \vert}				=
\frac{\sqrt{k^2 +1}} {k},
\end{align*}
which finishes the proof.
\end{proof} 

With the same technique and somewhat more effort one can give an embedding $\gamma_q$ of $T_{2,q}$ with $\delta(\gamma_q) \sim \frac{\pi}{2} \frac{q}{\log q}$. In addition a more technical proof yields that this asymptotical upper bound for $\delta(T_{2,q})$ is sharp for those embeddings of $T_{2,q}$ that project to a standard knot diagram via a linear projection. This let the author to the following.
\begin{question}
Is $\delta(T_{2,q})$ up to a constant asymptotically equal to $q/ \log q$? And if yes, is the constant equal to $\pi/2$?
\end{question}



\end{document}